\documentclass[12pt]{article}%
\usepackage{amssymb}
\usepackage{amsfonts}
\usepackage{graphicx}
\usepackage{amsmath}
\usepackage{hyperref}%
\providecommand{\U}[1]{\protect\rule{.1in}{.1in}}
\newtheorem{theorem}{Theorem}

\newtheorem{proposition}[theorem]{Proposition}
\newtheorem{remark}[theorem]{Remark}

\newenvironment{proof}[1][Proof]{\textbf{#1.} }{\ \rule{0.5em}{0.5em}}
\begin{document}

\title{A new approach to Tur\'{a}n inequalities for classical orthogonal polynomials}
\author{O. Kounchev and H. Render}
\maketitle

\begin{abstract}
Let $f,g\in C^{2}\left(  a,b\right)  .$ In this paper we investigate
properties of the function%
\[
v_{g}\left(  f\right)  \left(  x\right)  =-\frac{f^{2}\left(  x\right)
}{g\left(  x\right)  }\frac{d}{dx}\frac{g\left(  x\right)  f^{\prime}\left(
x\right)  }{f\left(  x\right)  }%
\]
where $f:\left(  a,b\right)  \rightarrow\mathbb{R}$ is a solution of the
second order differential equation $f^{\prime\prime}+pf^{\prime}+qf=0.$
Criteria for the monotonicity of $v_{g}\left(  f\right)  $ are given which
depend on the coefficients $p$ and $q$. It is shown that for classical
orthogonal polynomials $p_{n}\left(  x\right)  $ the Tur\'{a}n inequalities
are equivalent to the non-negativity of the functions $v_{g}\left(
p_{n}\right)  $ for a suitable function $g.$ This allows new and simple proofs
of classical and recent Tur\'{a}n inequalities.

\end{abstract}



{\large CORRESPONDING AUTHOR:} Ognyan Kounchev,

email:\ kounchev@abv.bg;

\textbf{Mailing address}: Institute of Mathematics and Informatics, Bulgarian
Academy of Sciences, Acad. G. Bonchev Str., bl. 8, 1113 Sofia, Bulgaria

{\large AUTHOR:} Hermann Render, email:\ hermann.render@ucd.ie;

\textbf{Mailing address}: School of Mathematics and Statistics, University
College Dublin, Belfield 4, Dublin, Ireland.

\bigskip

\bigskip

\section{Introduction}

Log-concavity and log-convexity of functions play an important role in many
areas of mathematics, in particular for the analysis of distribution functions
in statistics. Let us recall that a \emph{positive} function $f:\left(
a,b\right)  \rightarrow\left(  0,\infty\right)  $ is called log-concave
(log-convex resp.) if the function $x\longmapsto\log f\left(  x\right)  $ is
concave (convex resp.), see also \cite{AVV06} and \cite{AVV07} for a deeper
discussion of various concepts of convexity and concavity. If $f $ is twice
differentiable then $f$ is log-concave if and only if
\[
\frac{d^{2}}{dx^{2}}\log f\left(  x\right)  =\frac{d}{dx}\frac{f^{\prime
}\left(  x\right)  }{f\left(  x\right)  }=\frac{f^{\prime\prime}\left(
x\right)  f\left(  x\right)  -f^{\prime}\left(  x\right)  f^{\prime}\left(
x\right)  }{f\left(  x\right)  ^{2}}\leq0
\]
for all $x\in\left(  a,b\right)  $. This is equivalent to the statement that
$f^{\prime}\left(  x\right)  /f\left(  x\right)  $ is decreasing. From a
mathematical point it has advantages to investigate positivity and
monotonicity properties of the function
\begin{equation}
v\left(  f\right)  :=\left(  f^{\prime}\right)  ^{2}-f^{\prime\prime}f
\label{eqdefv}%
\end{equation}
which is meaningful for any twice differentiable \emph{real-valued} $f:\left(
a,b\right)  \rightarrow\mathbb{R}$. Note that $v\left(  f\right)  \left(
x\right)  \geq0$ implies that the function $\log\left\vert f\left(  x\right)
\right\vert $ is log-concave and $x\longmapsto\frac{f^{\prime}\left(
x\right)  }{f\left(  x\right)  }$ is decreasing on each subinterval which does
not contain a zero of $f$. Our methods of proof allow an important
generalization of the function $v\left(  f\right)  $: for a function $g\in
C^{2}\left(  a,b\right)  $ define%

\begin{equation}
v_{g}\left(  f\right)  \left(  x\right)  :=-\frac{f\left(  x\right)  ^{2}%
}{g\left(  x\right)  }\frac{d}{dx}\frac{g\left(  x\right)  f^{\prime}\left(
x\right)  }{f\left(  x\right)  }.\label{eqdefvg}%
\end{equation}
It is easy to see that $v_{g}\left(  f\right)  \left(  x\right)  $ is well
defined for all $x\in\left(  a,b\right)  $ for which $g\left(  x\right)
\neq0$ and $f\left(  x\right)  \neq0.$ Note that $v_{g}\left(  f\right)
=v\left(  f\right)  $ if $g$ is the constant function $1.$ When $v_{g}\left(
f\right)  \left(  x\right)  \geq0$ and $g\left(  x\right)  >0$ for all
$x\in\left(  a,b\right)  $, then $x\longmapsto\frac{g\left(  x\right)
f^{\prime}\left(  x\right)  }{f\left(  x\right)  }$is decreasing on each
subinterval which does not contain zeros of $f.$ An important example is the
elasticity of a function $f:\left(  a,b\right)  \rightarrow\left(
0,\infty\right)  $ defined by
\[
\text{El}\left(  f\right)  \left(  x\right)  =\frac{xf^{\prime}\left(
x\right)  }{f\left(  x\right)  }%
\]
which is an important tool  in econometrics. Clearly the elasticity of a
function $f:\left(  0,b\right)  \rightarrow\left(  0,\infty\right)  $ is
decreasing if and only if $v_{g}\left(  f\right)  \left(  x\right)  \geq0$ for
all $x\in\left(  0,b\right)  $ with $g\left(  x\right)  =x.$

An elementary but important fact is the following multiplicativity property of
$v_{g}$: If $v_{g}\left(  f_{j}\right)  \geq0$ on $\left(  a,b\right)  $ for
$j=1,...,n$ then $v_{g}\left(  f\right)  \geq0$ on $\left(  a,b\right)  $ for
the product $f=f_{1}\cdots f_{n}.$ Applied to the case that $f_{j}\left(
x\right)  =x-a_{j}$ are linear polynomials for $j=1,...,n,$ and $g\left(
x\right)  =1$ (so $v_{g}=v)$ it follows that $v\left(  P_{n}\right)  \left(
x\right)  \geq0$ for all $x\in\mathbb{R}$ whenever $P_{n}\left(  x\right)  $
is a polynomial with only real zeros since%

\[
v\left(  f_{j}\right)  \left(  x\right)  =-\left(  x-a_{j}\right)  ^{2}%
\frac{d}{dx}\frac{1}{x-a_{j}}=1.
\]

In the beginning of the 1940's, in a seminal paper, Tur\'{a}n proved that for
the Legendre polynomial $P_{n}\left(  x\right)  $ of degree $n$ the inequality%
\[
\Delta\left(  P_{n}\right)  \left(  x\right)  :=P_{n}\left(  x\right)
^{2}-P_{n+1}\left(  x\right)  P_{n-1}\left(  x\right)  \geq0
\]
holds for all $x\in\left[  -1,1\right]  ,$ see \cite{Tura50}. G. Szeg\"{o} has
given four proofs in \cite{Szeg48}. One proof was based on results about the
Laguerre-P\'{o}lya class and allowed an extension of the Tur\'{a}n inequality
to Hermite, Laguerre, ultrapsherical polynomials. Important contributions to
Tur\'{a}n inequalities were then obtained by several authors, see e.g.
\cite{Szas50}, \cite{ThNa51}, \cite{VeRa57}.

Although there is strong similarity between $v\left(  P_{n}\right)  $ and
$\Delta\left(  P_{n}\right)  $ (noting that lowering the index $n$ is
comparable to differentiation) the inequalities $v\left(  P_{n}\right)  \geq0$
and $\Delta\left(  P_{n}\right)  \geq0$ are not comparable: for Legendre
polynomials the error function $v\left(  P_{n}\right)  -\Delta\left(
P_{n}\right)  $ is oscillating. Nonetheless, in this paper we show that there
is an intimate relationship: for the normalized ultraspherical polynomials
$p_{n}^{\left(  \lambda\right)  }\left(  x\right)  $ we prove the new and
remarkable identity
\[
\Delta\left(  p_{n}^{\left(  \lambda\right)  }\right)  =\frac{\left(
1-x^{2}\right)  ^{2}}{n\left(  n+2\lambda\right)  }v_{g}\left(  p_{n}^{\left(
\lambda\right)  }\right)
\]
where $g=\left(  1-x^{2}\right)  ^{\frac{1}{2}}.$ Using the multiplicativity
property of $v_{g}$ the Tur\'{a}n inequality for ultraspherical polynomials
can be established just by verifying that for each linear polynomial
$f_{j}\left(  x\right)  =x-a_{j}$ with $a_{j}\in\left[  -1,1\right]  $ the expression%

\[
v_{g}\left(  f_{j}\right)  \left(  x\right)  =-\frac{\left(  x-a_{j}\right)
^{2}}{\left(  1-x^{2}\right)  ^{\frac{1}{2}}}\frac{d}{dx}\frac{\left(
1-x^{2}\right)  ^{\frac{1}{2}}}{x-a_{j}}=\frac{1-a_{j}x}{1-x^{2}}\geq0
\]
holds for all $x\in\left[  -1,1\right]  $. In the same spirit we shall give a
new and shorter proof of an interesting Tur\'{a}n inequality of G. Nikolov and
V. Pillwein \cite{NiPi15} stating that
\[
\left\vert x\right\vert \left(  p_{n}^{\left(  \lambda\right)  }\right)
^{2}-p_{n+1}^{\left(  \lambda\right)  }\left(  x\right)  p_{n-1}^{\left(
\lambda\right)  }\left(  x\right)  \geq0
\]
for all $x\in\left[  -1,1\right]  ,$ and $-\frac{1}{2}<\lambda\leq\frac{1}%
{2}.$ For $\lambda=\frac{1}{2},$ this inequality was earlier proved by Gerhold
and Kauers in \cite{GeKa06} for the case of Legendre polynomials.  We shall
also discuss extensions of this result due to I. Krasikov \cite{Kras23}.

In this paper the following general question is addressed: how can we prove
positivity of the function $v_{g}\left(  f\right)  $ defined in (\ref{eqdefvg}%
) on the interval $\left(  a,b\right)  $ when we just assume that $f\in
C^{2}\left(  a,b\right)  ?$ An elegant way is to find a suitable function
$H\left(  x\right)  $ such that the product $H\left(  x\right)  v_{g}\left(
f\right)  \left(  x\right)  $ is an increasing or decreasing function. Then
one can use limit arguments at the end points of the interval $\left(
a,b\right)  $ to infer positivity or negativity of the function $H\left(
x\right)  v_{g}\left(  f\right)  \left(  x\right)  $ which leads to positivity
results of the function $v_{g}\left(  f\right)  \left(  x\right)  $. Finding a
suitable function $H\left(  x\right)  $ seems to require non-trivial
heuristics arguments, and one needs more special information about the
function $f.$ In this paper we shall assume that $f$ is a solution of a second
order differential equation
\begin{equation}
f^{\prime\prime}+pf^{\prime}+qf=0\label{eqODE}%
\end{equation}
with $p,q\in C^{2}\left(  a,b\right)  .$ Replacing $f^{\prime\prime
}=-pf^{\prime}-qf$ in the definition (\ref{eqdefv}) leads to the simple
identity
\begin{equation}
v\left(  f\right)  =f^{\prime}f^{\prime}+pf^{\prime}f+qf^{2}.\label{eqvf}%
\end{equation}
For $v_{g}\left(  f\right)  $ we shall obtain a similar formula: an
application of the product rule shows that
\begin{equation}
v_{g}\left(  f\right)  =-\frac{1}{g}f^{2}\cdot\left(  g^{\prime}%
\frac{f^{\prime}}{f}+g\frac{d}{dx}\frac{f^{\prime}}{f}\right)  =v\left(
f\right)  -\frac{g^{\prime}}{g}f^{\prime}f.\label{eqvfalpha}%
\end{equation}
Using (\ref{eqvfalpha}) and (\ref{eqvf}) one obtains the identity
\[
v_{g}\left(  f\right)  =f^{\prime}f^{\prime}+\left(  p-\frac{g^{\prime}}%
{g}\right)  f^{\prime}f+qf^{2}.
\]
Now define $P=\int p\left(  x\right)  dx$. In order to prove positivity of
$v_{g}\left(  f\right)  $ we define the factor
\[
H\left(  x\right)  =\frac{e^{P\left(  x\right)  }}{p_{g}\left(  x\right)
}\text{ with }p_{g}\left(  x\right)  =p\left(  x\right)  -\frac{g^{\prime
}\left(  x\right)  }{g\left(  x\right)  }.
\]
This choice of $H\left(  x\right)  $ is based on the following identity proven
in Section 3 of the present paper:
\begin{equation}
p_{g}^{2}e^{-P}\frac{d}{dx}\left(  \frac{e^{P}}{p_{g}}v_{g}\left(  f\right)
\right)  =-\left(  \frac{g^{\prime}}{g}p_{g}+p_{g}^{\prime}\right)  f^{\prime
}f^{\prime}-\left(  qp_{g}^{\prime}-q\frac{g^{\prime}}{g}p_{g}-q^{\prime}%
p_{g}\right)  f^{2}.\label{eqMainid}%
\end{equation}
This rather technical formula simplifies when we take $g=1.$ Then $p_{g}%
=p\neq0$, and it follows from (\ref{eqMainid}) that $\frac{e^{P}}{p}v\left(
f\right)  $ is increasing if
\begin{equation}
\text{ }p^{\prime}\leq0\text{ and }p^{\prime}q-q^{\prime}p\leq0\text{ on
}\left(  a,b\right)  .\label{eqcritb}%
\end{equation}
If $\frac{e^{P}}{p}v\left(  f\right)  \left(  x\right)  \rightarrow0$ for
$x\rightarrow a$ then one obtains positivity of $\frac{e^{P}}{p}v\left(
f\right)  $ on $\left(  a,b\right)  .$ Reversing the inequalities
(\ref{eqcritb}) leads to a sufficient condition that $\frac{e^{P}}{p}v\left(
f\right)  $ is decreasing.

The paper is organized as follows: in Section $2$  we present basic properties
of the function $v_{g}\left(  f\right)  .$ In Section 3 we will prove
identity (\ref{eqMainid}). In Section 4 we discuss  Tur\'{a}n inequality for
the Hermite polynomial $H_{n}$ and remark (as is known for a long time) that
\[
\Delta\left(  H_{n}\right)  =H_{n}^{2}-H_{n-1}H_{n+1}=\frac{1}{2n}v\left(
H_{n}\right)  .
\]
Hermite polynomials are special cases of Hermite functions $H_{\nu}$ (see
Section 10.2 in \cite{Lebe} for a nice exposition) which are entire functions
defined for all complex parameters $\nu.$ We use the identity (\ref{eqMainid})
to show that for all positive $\nu$ (in particular for natural numbers $n)$
the inequality
\[
v\left(  H_{\nu}\right)  \left(  x\right)  >0
\]
holds for all $x\in\left[  -1,1\right]  .$ For $\nu<0$ it is known that
$H_{\nu}\left(  x\right)  >0$ for all real $x$ and we shall show that
$v\left(  H_{\nu}\right)  \left(  x\right)  <0$ for all real $x,$ i.e. that
$H_{\nu}$ is a log-convex which was proved in \cite{Koch20} using techniques
from stochastic processes. Section 5 is devoted to generalized Laguerre
polynomials. In Section 6 and 7 we discuss  Tur\'{a}n inequality for
ultraspherical polynomials. Finally we mention that our methods can also be
applied to Jacobi polynomials which we plan to discuss in a future paper.

Let us emphasize that there is a vast amount of literature about Tur\'{a}n
inequalities for special functions, see e.g. \cite{Bari06}, \cite{BPV11}. In a
forth-coming paper we shall show that the identity (\ref{eqMainid}) allows us
to give a systematic and unified approach for proving log-concavity and
Tur\'{a}n inequalities for Bessel functions $J_{\nu}$ of the first kind, and
modified Bessel functions $I_{\nu}$ of the first kind, see \cite[in
preparation]{KoRe27}.

\section{Basic facts}

Throughout the paper we assume that $f:\left(  a,b\right)  \rightarrow
\mathbb{R}$ is a twice differentiable real-valued function and the notation%
\[
v\left(  f\right)  =f^{\prime}f^{\prime}-f^{\prime\prime}f
\]
is used. Often we simply write $v$ instead of $v\left(  f\right)  $ if there
is no ambiguity. Further $v^{\prime}\left(  f\right)  $ will denote the
derivative of $v\left(  f\right)  .$ It is instructive to look at some
examples: it is easy to see that $v\left(  f\right)  =0$ for a constant
function and the exponential function $e^{x}$. Moreover, for $\alpha
\in\mathbb{R}$ and for $f_{\alpha}\left(  x\right)  =x^{\alpha}$
\begin{equation}
v\left(  f_{\alpha}\right)  \left(  x\right)  =\alpha^{2}x^{2\alpha-2}%
-\alpha\left(  \alpha-1\right)  x^{2\alpha-2}=\alpha x^{2\alpha-2}%
.\label{eqnupower}%
\end{equation}
Hence $v\left(  f_{\alpha}\right)  \left(  x\right)  \geq0$ for all
$\alpha\geq0$ and $x>0.$

\begin{proposition}
\label{Prop0}Let $g,f_{1},f_{2}\in C^{2}\left(  a,b\right)  .$ Then
\[
v_{g}\left(  f_{1}f_{2}\right)  =f_{2}^{2}v_{g}\left(  f_{1}\right)
+f_{1}^{2}v_{g}\left(  f_{2}\right)  .
\]

\end{proposition}

\begin{proof}
This follows by a straightforward calculation using definition (\ref{eqdefvg}%
):%
\begin{align*}
v_{g}\left(  f_{1}f_{2}\right)   &  =-\frac{1}{g}f_{1}^{2}f_{2}^{2}\frac
{d}{dx}\frac{g\frac{d}{dx}\left(  f_{1}f_{2}\right)  }{f_{1}f_{2}}=-\frac
{1}{g}f_{1}^{2}f_{2}^{2}\frac{d}{dx}\frac{g\left(  f_{1}^{\prime}f_{2}%
+f_{1}f_{2}^{\prime}\right)  }{f_{1}f_{2}}\\
&  =-\frac{1}{g}f_{1}^{2}f_{2}^{2}\frac{d}{dx}\left(  \frac{gf_{1}^{\prime}%
}{f_{1}}+\frac{gf_{2}^{\prime}}{f_{2}}\right)  =f_{2}^{2}v_{g}\left(
f_{1}\right)  +f_{1}^{2}v_{g}\left(  f_{2}\right)  .
\end{align*}

\end{proof}

Proposition \ref{Prop0} shows that $v_{g}\left(  f_{1}f_{2}\right)  \geq0$ on
$\left(  a,b\right)  $ if $v_{g}\left(  f_{1}\right)  \geq0$ and $v_{g}\left(
f_{2}\right)  \geq0$ on $\left(  a,b\right)  .$ Using induction one obtains:

\begin{theorem}
\label{ThmProduct}Let $g$ and $f_{j}\in C^{2}\left(  a,b\right)  $ for
$j=1,...,n.$ If $v_{g}\left(  f_{j}\right)  \geq0$ on $\left(  a,b\right)  $
and $j=1,...,n$ then $v_{g}\left(  f\right)  \geq0$ on $\left(  a,b\right)  $
for the product $f=f_{1}\cdots f_{n}.$
\end{theorem}

Since $v\left(  e^{x}\right)  =0$, Proposition \ref{Prop0} implies that
\[
v\left(  e^{x}f\right)  =e^{2x}v\left(  f\right)  +f^{2}v\left(  e^{x}\right)
=e^{2x}v\left(  f\right)  .
\]
Another nice consequence of Proposition \ref{Prop0} is the following
well-known result (credited to Laguerre in \cite{Cham17}) which usually is
proved by considering the partial fraction decomposition of $p^{\prime}/p$..

\begin{theorem}
\label{Thmpolynomial}Assume that $p\left(  x\right)  $ is a polynomial with
only real zeros. Then $v\left(  p\right)  \left(  x\right)  \geq0$ for all
real $x.$
\end{theorem}

\begin{proof}
Since $p\left(  x\right)  $ is a product of linear factors it suffices to show
$v\left(  f\right)  \geq0$ for $f\left(  x\right)  =A\left(  x-a\right)  $
with $a,A$ real. Clearly $v\left(  f\right)  =A^{2}\geq0$.
\end{proof}

It is easy to see that $v\left(  p\right)  \left(  x_{0}\right)  =0$ if
$x_{0}$ is a double zero of $p$, so one can not prove positivity in Theorem
\ref{Thmpolynomial}.

\begin{theorem}
Assume that $p\left(  x\right)  $ is a polynomial of degree $\geq1$ with only
real and distinct zeros. Then $v\left(  p\right)  \left(  x\right)  >0$ for
all real $x.$
\end{theorem}

\begin{proof}
The proof follows by induction over the degree $n\geq1$: For $n=1$ let
$f\left(  x\right)  =ax+b$ with $a\neq0.$ Then $v\left(  f\right)  =a^{2}>0.$
If $p_{n+1}$ is a polynomial of degree $n+1$ with only real and distinct,
write it as $p_{n+1}\left(  x\right)  =\left(  x-a\right)  p_{n}\left(
x\right)  $ where $a$ is a real number distinct from the real and distinct
zeros of $p_{n}.$ By induction hypothesis $v\left(  p_{n}\right)  $ is
positive, and $v\left(  p_{n+1}\right)  \left(  x\right)  =\left(  x-a\right)
^{2}v\left(  p_{n}\right)  \left(  x\right)  +p_{n}^{2}\left(  x\right)  $
$>0$ for all $x\neq a.$ If $x=a$ we see that $p_{n}^{2}\left(  a\right)  >0$
since $a$ is distinct from the zeros of $p_{n}.$
\end{proof}

Let us mention a connection to the Laguerre-P\'{o}lya class $\mathcal{L}%
-\mathcal{P}$, the class of all entire functions with real coefficients which
can be approximated by real polynomials having only real zeros. A basic fact
states that a function $f\in\mathcal{L}-\mathcal{P}$ satisfies the
\emph{Laguerre inequality}
\[
v\left(  f\right)  \left(  x\right)  =f^{\prime}\left(  x\right)  f^{\prime
}\left(  x\right)  -f\left(  x\right)  f^{\prime\prime}\left(  x\right)
\geq0\text{ for all }x\in\mathbb{R}.
\]
Since the derivative $f^{\prime}$ of a function in $\mathcal{L}-\mathcal{P}$
is again in $\mathcal{L}-\mathcal{P}$ one obtains $v\left(  f^{\left(
k\right)  }\right)  \geq0$ for the $k$-th derivative $f^{\left(  k\right)  }$.
We refer to \cite{Patr73} and \cite{CsEs05} for these facts and further
interesting higher order Tur\'{a}n inequalities.

The following result is straightforward consequence of the chain rule and the
proof is omitted.

\begin{proposition}
\label{PropEven} Let $f\in C^{2}\left(  \left(  -a,a\right)  \right)  $ be
even or odd function. Then $v\left(  f\right)  $ is even.
\end{proposition}

\section{Results for solutions of $f^{\prime\prime}+pf^{\prime}+qf=0$}

The following result can be found in \cite{KoRe26}. For convenience of the
reader we include a short proof:

\begin{proposition}
\label{Prop1}Let $p,q\in C^{1}\left(  a,b\right)  .$ If $f$ is a solution of
$f^{\prime\prime}+pf^{\prime}+qf=0$ then%
\begin{equation}
v^{\prime}\left(  f\right)  +pv\left(  f\right)  =p^{\prime}f^{\prime
}f+q^{\prime}f^{2}. \label{eqvde}%
\end{equation}

\end{proposition}

\begin{proof}
Substitute $f^{\prime\prime}=-pf^{\prime}-qf$ in the definition of $v=\left(
f^{\prime}\right)  ^{2}-f^{\prime\prime}f,$ and one obtains
\begin{equation}
v=\left(  f^{\prime}\right)  ^{2}+pf^{\prime}f+qf^{2}. \label{eq1}%
\end{equation}
Differentiation of (\ref{eq1}) gives
\begin{equation}
v^{\prime}=2f^{\prime}f^{\prime\prime}+p^{\prime}f^{\prime}f+pf^{\prime\prime
}f+pf^{\prime}f^{\prime}+q^{\prime}f^{2}+2qf^{\prime}f. \label{eq2}%
\end{equation}
Replace $f^{\prime\prime}$ by $-pf^{\prime}-qf$, so
\begin{align*}
v^{\prime}  &  =-2pf^{\prime}f^{\prime}-2qf^{\prime}f+p^{\prime}f^{\prime
}f-p^{2}f^{\prime}f-pqf^{2}+pf^{\prime}f^{\prime}+q^{\prime}f^{2}+2qf^{\prime
}f\\
&  =-pf^{\prime}f^{\prime}+\left(  p^{\prime}-p^{2}\right)  f^{\prime
}f+\left(  q^{\prime}-pq\right)  f^{2}.
\end{align*}
Using (\ref{eq1}) we obtain our claim that $v^{\prime}=-pv+p^{\prime}%
f^{\prime}f+q^{\prime}f^{2}.$
\end{proof}

In the introduction we have shown that%
\begin{equation}
v_{g}\left(  f\right)  =v\left(  f\right)  -\frac{g^{\prime}}{g}f^{\prime}f.
\label{eqvalpha}%
\end{equation}
Since $v\left(  f\right)  =\left(  f^{\prime}\right)  ^{2}+pf^{\prime}%
f+qf^{2}$ it follows that
\begin{equation}
v_{g}\left(  f\right)  =f^{\prime}f^{\prime}+p_{g}f^{\prime}f+qf^{2}\text{
where }p_{g}:=p-\frac{g^{\prime}}{g} \label{eqAlpha}%
\end{equation}

\begin{proposition}
\label{Prop1alpha}Let $p,q\neq0\in C^{1}\left(  a,b\right)  $ and $g\in
C^{2}\left(  a,b\right)  .$ Define $p_{g}=p-\frac{g^{\prime}}{g}$. If $f$ is a
solution of $f^{\prime\prime}+pf^{\prime}+qf=0$ then%
\begin{equation}
v_{g}^{\prime}\left(  f\right)  +pv_{g}\left(  f\right)  =-\frac{g^{\prime}%
}{g}f^{\prime}f^{\prime}+p_{g}^{\prime}f^{\prime}f+\left(  q^{\prime}%
+q\frac{g^{\prime}}{g}\right)  f^{2} \label{eqvdealpha}%
\end{equation}

\end{proposition}

\begin{proof}
We know $v_{g}=v-\frac{g^{\prime}}{g}f^{\prime}f.$ It follows that%
\begin{equation}
v_{g}^{\prime}+pv_{g}=v^{\prime}-\frac{d}{dx}\left(  \frac{g^{\prime}}%
{g}f^{\prime}f\right)  +pv-p\frac{g^{\prime}}{g}f^{\prime}f \label{eqvg}%
\end{equation}
By Proposition \ref{Prop1} we know that $v^{\prime}+pv=p^{\prime}f^{\prime
}f+q^{\prime}f^{2}$. Further
\[
-\frac{d}{dx}\left(  \frac{g^{\prime}}{g}f^{\prime}f\right)  =-\frac{d}%
{dx}\left(  \frac{g^{\prime}}{g}\right)  f^{\prime}f-\frac{g^{\prime}}%
{g}\left(  f^{\prime\prime}f+f^{\prime}f^{\prime}\right)  .
\]
Replace $f^{\prime\prime}=-pf^{\prime}-qf$ and note that the last term in
(\ref{eqvg}), namely $-p\frac{g^{\prime}}{g}f^{\prime}f$, cancels out leading
to
\begin{align*}
v_{g}^{\prime}+pv_{g}  &  =p^{\prime}f^{\prime}f+q^{\prime}f^{2}-\frac{d}%
{dx}\left(  \frac{g^{\prime}}{g}\right)  f^{\prime}f+\frac{g^{\prime}}%
{g}qff-\frac{g^{\prime}}{g}f^{\prime}f^{\prime}\\
&  =-\frac{g^{\prime}}{g}f^{\prime}f^{\prime}+\left(  p^{\prime}-\frac{d}%
{dx}\left(  \frac{g^{\prime}}{g}\right)  \right)  f^{\prime}f+\left(
q^{\prime}+q\frac{g^{\prime}}{g}\right)  f^{2}.
\end{align*}
Since $p_{g}=p-\frac{g^{\prime}}{g}$ we see that the coefficient of
$f^{\prime}f$ is equal to $p_{g}^{\prime}.$
\end{proof}

The following is a central result in our paper:

\begin{theorem}
\label{ThmMainODE}Let $p,q\in C^{1}\left(  a,b\right)  $ and define $P=\int
p\left(  x\right)  dx$. For $g\in C^{2}\left(  a,b\right)  $ define
$p_{g}=p-\frac{g^{\prime}}{g}.$ If $f$ is a solution of $f^{\prime\prime
}+pf^{\prime}+qf=0$ and $p_{g}\neq0$ then the identity
\[
p_{g}^{2}e^{-P}\frac{d}{dx}\left(  \frac{e^{P}}{p_{g}}v_{g}\left(  f\right)
\right)  =-\left(  p_{g}^{\prime}+\frac{g^{\prime}}{g}p_{g}\right)  f^{\prime
}f^{\prime}-\left(  qp_{g}^{\prime}-q\frac{g^{\prime}}{g}p_{g}-q^{\prime}%
p_{g}\right)  f^{2}.
\]
holds. If $p_{g}=0$ then
\[
e^{-2P}\frac{d}{dx}\left(  e^{2P}v_{g}\left(  f\right)  \right)  =\left(
q^{\prime}+2pq\right)  f^{2}.
\]

\end{theorem}

\begin{proof}
Multiply $v_{g}=f^{\prime}f^{\prime}+p_{g}f^{\prime}f+qf^{2}$ with
$p_{g}^{\prime},$ so we obtain
\[
p_{g}^{\prime}v_{g}=p_{g}^{\prime}f^{\prime}f^{\prime}+p_{g}p_{g}^{\prime
}f^{\prime}f+qp_{g}^{\prime}f^{2}.
\]
Replace $p_{g}^{\prime}f^{\prime}f$ in the middle term by using the equation
(\ref{eqvdealpha}), namely
\[
v_{g}^{\prime}+pv_{g}=-\frac{g^{\prime}}{g}f^{\prime}f^{\prime}+p_{g}^{\prime
}f^{\prime}f+\left(  q^{\prime}+q\frac{g^{\prime}}{g}\right)  f^{2}.
\]
Then
\begin{align*}
p_{g}^{\prime}v_{g}  &  =p_{g}^{\prime}f^{\prime}f^{\prime}+p_{g}\left(
v_{g}^{\prime}+pv_{g}+\frac{g^{\prime}}{g}f^{\prime}f^{\prime}-\left(
q^{\prime}+q\frac{g^{\prime}}{g}\right)  f^{2}\right)  +qp_{g}^{\prime}f^{2}\\
&  =p_{g}^{\prime}f^{\prime}f^{\prime}+p_{g}v_{g}^{\prime}+p_{g}pv_{g}%
+\frac{g^{\prime}}{g}p_{g}f^{\prime}f^{\prime}-p_{g}\left(  q^{\prime}%
+q\frac{g^{\prime}}{g}\right)  f^{2}+qp_{g}^{\prime}f^{2}.
\end{align*}
It follows that
\[
p_{g}v_{g}^{\prime}+p_{g}pv_{g}-p_{g}^{\prime}v_{g}=-\left(  p_{g}^{\prime
}+\frac{g^{\prime}}{g}p_{g}\right)  f^{\prime}f^{\prime}-\left(
qp_{g}^{\prime}-qp_{g}\frac{g^{\prime}}{g}-q^{\prime}p_{g}\right)  f^{2}.
\]
Let $r$ be a function and $R$ an anti-derivative of $r.$ Then
\[
e^{-R}\frac{d}{dx}\left(  e^{R}F\right)  =e^{-R}\left(  e^{R}rF+e^{R}%
F^{\prime}\right)  =F^{\prime}+rF.
\]
Put $F=v_{g}$ and $r=p-\frac{p_{g}^{\prime}}{p_{g}}$ and $\int rdx=P-\ln
p_{g}$ so we have $e^{R}=\frac{e^{P}}{p_{g}}.$ Then
\[
p_{g}e^{-P}\frac{d}{dx}\left(  \frac{e^{P}}{p_{g}}v_{g}\right)  =v_{g}%
^{\prime}+pv_{g}-\frac{p_{g}^{\prime}}{p_{g}}v_{g}.
\]
The first statement of the theorem follows.

2. In the case $p_{g}=0$ we have $v_{g}=f^{\prime}f^{\prime}+qf^{2}$ by
(\ref{eqAlpha}). Then
\[
v_{g}^{\prime}=2f^{\prime\prime}f^{\prime}+q^{\prime}f^{2}+2qf^{\prime
}f=2\left(  -pf^{\prime}-qf\right)  f^{\prime}+q^{\prime}f^{2}+2qf^{\prime
}f=-2pf^{\prime}f^{\prime}+q^{\prime}f^{2}.
\]
It follows that
\[
v_{g}^{\prime}+2pv_{g}=-2pf^{\prime}f^{\prime}+q^{\prime}f^{2}+2pf^{\prime
}f^{\prime}+2pqf^{2}=\left(  q^{\prime}+2pq\right)  f^{2}.
\]
Hence
\[
e^{-2P}\frac{d}{dx}\left(  e^{2P}v_{g}\right)  =v_{g}^{\prime}+2pv_{g}=\left(
q^{\prime}+2pq\right)  f^{2}.
\]

\end{proof}

It is instructive to consider the case $g=1,$ so $\frac{g^{\prime}}{g}=0$ and
$p_{g}=p.$ Then
\begin{equation}
pe^{-P}\frac{d}{dx}\left(  \frac{e^{P}}{p}v\right)  =-\frac{p^{\prime}}%
{p}f^{\prime}f^{\prime}-\left(  q\frac{p^{\prime}}{p}-q^{\prime}\right)
f^{2}. \label{eqgeq1}%
\end{equation}
Thus the function $\frac{e^{P}}{p}v$ is increasing (decreasing resp.) if
$\frac{p^{\prime}}{p}$ and $q\frac{p^{\prime}}{p}-q^{\prime}$ are both
negative (positive resp.). Thus we have

\begin{theorem}
\label{Thm2Monotdiffeq}Let $p,q\in C^{1}\left(  a,b\right)  $ such that
$p\left(  x\right)  \neq0$ for all $x\in\left(  a,b\right)  $ and define
$P=\int p\left(  x\right)  dx$. If $f\in C^{2}\left(  a,b\right)  $ is a
solution of $f^{\prime\prime}+pf^{\prime}+qf=0,$ and if%
\[
\text{ }p^{\prime}\leq0\text{ and }p^{\prime}q-q^{\prime}p\leq0\text{ on
}\left(  a,b\right)
\]
then $\frac{e^{P}}{p}v\left(  f\right)  $ is increasing on $\left(
a,b\right)  $, and if%
\[
\text{ }p^{\prime}\geq0\text{ and }p^{\prime}q-q^{\prime}p\geq0\text{ on
}\left(  a,b\right)
\]
then $\frac{e^{P}}{p}v\left(  f\right)  $ is decreasing on $\left(
a,b\right)  $.
\end{theorem}

\begin{remark}
If $p^{\prime}<0$ and $p^{\prime}q-q^{\prime}p<0$ on $\left(  a,b\right)  $
and $f\neq0$ then $\frac{e^{P}}{p}v\left(  f\right)  $ is strictly increasing
\end{remark}

\begin{proof}
Suppose that $\frac{d}{dx}\left(  \frac{e^{P}}{p}v\right)  \left(  x\right)
=0.$ Then both summands in (\ref{eqgeq1}) are non-negative. It follows that
$-p^{\prime}\left(  x\right)  f^{\prime}\left(  x\right)  ^{2}=0$ and
$-(p^{\prime}\left(  x\right)  q\left(  x\right)  -p\left(  x\right)
q^{\prime}\left(  x\right)  )f^{2}\left(  x\right)  =0.$ Our positivity
assumption implies $f\left(  x\right)  =f^{\prime}\left(  x\right)  =0.$ Since
$f$ is a solution of a second order differential equation, $f$ is identically
zero, a contradiction.
\end{proof}

\section{Hermite polynomials and Hermite functions}

Hermite polynomial $H_{n}\left(  x\right)  $ of degree $n$ satisfies the
differential equation%
\begin{equation}
f^{\prime\prime}-2xf^{\prime}+2nf=0.\label{eqODEHermite}%
\end{equation}
Since $H_{n}$ has only real and distinct zeros it follows from Theorem
\ref{Thmpolynomial} that
\[
v\left(  H_{n}\right)  \left(  x\right)  =H_{n}^{\prime}\left(  x\right)
^{2}-H_{n}^{\prime\prime}\left(  x\right)  H_{n}\left(  x\right)  >0
\]
for all real numbers $x$. Since $H_{n}$ is either even or odd it follows by
Proposition \ref{PropEven} that $v\left(  H_{n}\right)  $ is an even
polynomial of degree $2n-2.$ It was already observed in \cite{MuNa51} , cf.
\cite[p. 71]{Skov54}, that
\begin{equation}
\Delta\left(  H_{n}\right)  =H_{n}^{2}-H_{n-1}H_{n+1}=\frac{1}{2n}v\left(
H_{n}\right)  .\label{eqbyproduct}%
\end{equation}
Krasikov proved in \cite[Theorem 3]{Kras23} a sharper inequality which we
shall prove in a shorter way using our new techniques. As a byproduct we
obtain (\ref{eqbyproduct}).

\begin{theorem}
\label{ThmvHermite}Let $H_{n}$ be the Hermite polynomial of degree $n.$ Then
for all real $x,$ and
\[
\Delta\left(  H_{n}\right)  \left(  x\right)  \geq\frac{x^{2}}{x^{2}+\frac
{1}{2}}H_{n}^{2}\left(  x\right)  -H_{n-1}\left(  x\right)  H_{n+1}\left(
x\right)  \geq0.
\]

\end{theorem}

\begin{proof}
Since $H_{n}$ satisfies the differential equation (\ref{eqODEHermite}) we
infer by (\ref{eq1}) that%
\[
v\left(  H_{n}\right)  =H_{n}^{\prime}H_{n}^{\prime}-2xH_{n}^{\prime}%
H_{n}+2nH_{n}^{2}.
\]
The recurrence relations $H_{n}^{\prime}=2nH_{n-1}$ and $H_{n+1}%
=2xH_{n}-2nH_{n-1}$ show that%
\[
-H_{n-1}H_{n+1}=-\frac{1}{2n}H_{n}^{\prime}\left(  2xH_{n}-2nH_{n-1}\right)
=-\frac{2x}{2n}H_{n}^{\prime}H_{n}+H_{n}^{\prime}H_{n-1}.
\]
Apply again the relation $H_{n}^{\prime}=2nH_{n-1},$ and add $H_{n}^{2}.$ Then
we obtain (\ref{eqbyproduct}), namely
\[
\Delta\left(  H_{n}\right)  =\frac{1}{2n}H_{n}^{\prime}H_{n}^{\prime}-\frac
{x}{n}H_{n}^{\prime}H_{n}+H_{n}^{2}=\frac{1}{2n}v\left(  H_{n}\right)  .
\]
When we add $\frac{x^{2}}{x^{2}+\frac{1}{2}}H_{n}^{2}\left(  x\right)  $
instead of $H_{n}^{2}$ we obtain%
\[
D_{n}:=\frac{x^{2}}{x^{2}+\frac{1}{2}}H_{n}^{2}\left(  x\right)
-H_{n+1}\left(  x\right)  H_{n-1}\left(  x\right)  =\frac{x^{2}}{x^{2}%
+\frac{1}{2}}H_{n}^{2}-\frac{x}{n}H_{n}^{\prime}H_{n}+\frac{1}{2n}%
H_{n}^{\prime}H_{n}^{\prime}.
\]
Then
\begin{align*}
D_{n}  &  =\frac{1}{2n}\frac{x^{2}}{x^{2}+\frac{1}{2}}\left(  2nH_{n}%
^{2}-\frac{x^{2}+\frac{1}{2}}{x^{2}}2xH_{n}^{\prime}H_{n}+\frac{x^{2}+\frac
{1}{2}}{x^{2}}H_{n}^{\prime}H_{n}^{\prime}\right) \\
&  =\frac{1}{2n}\frac{x^{2}}{x^{2}+\frac{1}{2}}\left(  2nH_{n}^{2}%
-2xH_{n}^{\prime}H_{n}-\frac{1}{x}H_{n}^{\prime}H_{n}+H_{n}^{\prime}%
H_{n}^{\prime}+\frac{1}{2x^{2}}H_{n}^{\prime}H_{n}^{\prime}\right)  .
\end{align*}
Set now $g\left(  x\right)  =x,$ so $g^{\prime}\left(  x\right)  /g\left(
x\right)  =1/x.$ In view the identity $v_{g}\left(  f\right)  =f^{\prime
}f^{\prime}+\left(  p-\frac{g^{\prime}}{g}\right)  f^{\prime}f+qf^{2}$ we
obtain
\[
D_{n}=\frac{1}{2n}\frac{x^{2}}{x^{2}+\frac{1}{2}}\left(  v_{g}\left(
H_{n}\right)  +\frac{1}{2x^{2}}H_{n}^{\prime}H_{n}^{\prime}\right)  .
\]
Proposition \ref{PropEvenOdd} shows that $v_{g}\left(  H_{n}\right)  $ is
positive on the real line. Hence $D_{n}$ is positive for all $x.$
\end{proof}

The next result will be used later:

\begin{proposition}
\label{PropEvenOdd}Let $g\left(  x\right)  =x.$ Then $v_{g}\left(
H_{n}\right)  $ is positive on the real line.
\end{proposition}

\begin{proof}
Since $H_{n}$ is odd or even it is a product of factors $p_{a}\left(
x\right)  =A\left(  x-a\right)  \left(  x+a\right)  $ with real $a\neq0,$ and
in case of odd degree a factor $p\left(  x\right)  =x.$ Then
\[
v_{g}\left(  p_{2}\right)  \left(  x\right)  =-\frac{1}{x}\left(  A\left(
x-a\right)  \left(  x+a\right)  \right)  ^{2}\frac{d}{dx}\frac{x\frac{d}%
{dx}\left(  \left(  x-a\right)  \left(  x+a\right)  \right)  }{\left(
x-a\right)  \left(  x+a\right)  }=4a^{2}A^{2}>0
\]
and it is easy to see that $v_{g}\left(  p\right)  =0$ for $p\left(  x\right)
=x.$ By Theorem \ref{ThmProduct} $v_{g}\left(  H_{n}\right)  \left(  x\right)
\geq0$ for all real $x.$
\end{proof}

Hermite functions $H_{\nu}$ (see Section 10.2 in \cite{Lebe} for a nice
exposition) are entire functions defined for all complex parameters $\nu.$
They play an important role for problems using cylindrical coordinates in
physics (see \cite{Lebe}) but also appear in the context of numerical
analysis:
\[
H_{-1}\left(  x\right)  =\frac{\sqrt{\pi}}{2}e^{x^{2}}\left(
1-\operatorname{erf}\left(  x\right)  \right)
\]
where $\operatorname{erf}\left(  x\right)  =\frac{2}{\sqrt{\pi}}\int_{0}%
^{x}e^{-t^{2}}dt$ is the error function. The Hermite functions $H_{\nu}$ are
solutions of the differential equation%
\[
f^{\prime\prime}-2xf^{\prime}+2\nu f=0
\]
and they satisfy the relations $H_{\nu}^{\prime}=2\nu H_{\nu-1}$ and
$H_{\nu+1}=2xH_{\nu}-2\nu H_{\nu-1}$. Following the proof of Theorem
\ref{ThmvHermite} one has
\begin{equation}
\Delta\left(  H_{\nu}\right)  =\frac{1}{2\nu}H_{\nu}^{\prime}H_{\nu}^{\prime
}-\frac{x}{\nu}H_{\nu}^{\prime}H_{\nu}+H_{\nu}^{2}=\frac{1}{2\nu}v\left(
H_{\nu}\right)  .\label{eqDeltaH}%
\end{equation}
At first we prove

\begin{theorem}
Let $\nu>0.$ Then $v\left(  H_{\nu}\right)  \left(  x\right)  >0$ for all real
$x.$
\end{theorem}

\begin{proof}
We apply Theorem \ref{Thm2Monotdiffeq}: clearly $p^{\prime}=-2<0$ and
$p^{\prime}q-pq^{\prime}=-4\nu<0$ on the interval $\left(  -\infty
,\infty\right)  $. Further $P\left(  x\right)  =\int p\left(  x\right)
dx=-x^{2}$ and
\[
\frac{e^{P\left(  x\right)  }}{p\left(  x\right)  }=\frac{e^{-x^{2}}}{-2x}.
\]
By Theorem \ref{Thm2Monotdiffeq}, $\frac{e^{-x^{2}}}{-2x}v\left(  H_{\nu
}\right)  \left(  x\right)  $ is strictly increasing on $\left(
-\infty,0\right)  $ and $\left(  0,\infty\right)  .$ Thus for $x<0$ we have
\[
\lim_{x\rightarrow-\infty}\frac{e^{-x^{2}}}{-2x}v\left(  H_{\nu}\right)
\left(  x\right)  \leq\frac{e^{-x^{2}}}{-2x}v\left(  H_{\nu}\right)  \left(
x\right)
\]
and for $x>0,$we use that $\frac{e^{-x^{2}}}{2x}v\left(  H_{n}\right)  \left(
x\right)  $ is strictly decreasing on $\left(  0,\infty\right)  ,$ so
\[
\frac{e^{-x^{2}}}{2x}v\left(  H_{\nu}\right)  \left(  x\right)  \geq
\lim_{x\rightarrow\infty}\frac{e^{-x^{2}}}{2x}v\left(  H_{\nu}\right)  \left(
x\right)  .
\]
Using the asymptotic behaviour. of $H_{\nu}$, see equation (10.6.4) and
(10.6.7) in \cite{Lebe}, it follows that
\[
\lim_{x\rightarrow\pm\infty}\frac{e^{-x^{2}}}{2x}v\left(  H_{\nu}\right)
\left(  x\right)  =0.
\]
Then $\frac{e^{-x^{2}}}{-2x}v\left(  H_{\nu}\right)  \left(  x\right)  >0$ for
all $x<0$ and $\frac{e^{-x^{2}}}{2x}v\left(  H_{\nu}\right)  \left(  x\right)
>0$ for all $x>0.$ Thus $v\left(  H_{\nu}\right)  \left(  x\right)  >0$ for
all $x\neq0.$
\end{proof}

The following result was proved in \cite{Koch20} using techniques from
stochastic processes. We give an analytical proof in our framework.

\begin{theorem}
Let $\nu<0.$ Then $H_{\nu}\left(  x\right)  >0$ and $v\left(  H_{\nu}\right)
\left(  x\right)  <0$ for all real $x.$ Furthermore $\Delta\left(  H_{\nu
}\right)  \left(  x\right)  >0$ for all real $x.$
\end{theorem}

\begin{proof}
According to Proposition \ref{Prop1} we have $v^{\prime}\left(  f\right)
+pv\left(  f\right)  =p^{\prime}f^{\prime}f+q^{\prime}f^{2},$ hence with
$p=-2x$ and $q=2\nu$ and $f=H_{\nu}.$
\[
v^{\prime}\left(  f\right)  -2xv\left(  f\right)  =-2f^{\prime}f.
\]
Since $\nu<0$ it follows from the integral representation (10.5.7) in
\cite{Lebe} that $H_{\nu}\left(  x\right)  $ is a positive function for all
real $x,$and the same is true for $H_{\nu-1}.$ Since $H_{\nu}^{\prime}\left(
x\right)  =2\nu H_{\nu-1}\left(  x\right)  $ it follows $H_{\nu}^{\prime
}\left(  x\right)  $ is negative. Thus $H_{\nu}H_{\nu}^{\prime}$ is negative
and
\[
e^{x^{2}}\frac{d}{dx}\left(  e^{-x^{2}}v\left(  H_{\nu}\right)  \right)
=v^{\prime}\left(  H_{\nu}\right)  -2xv\left(  H_{\nu}\right)  =-2H_{\nu
}\left(  x\right)  H_{\nu}^{\prime}\left(  x\right)  >0
\]
for all real $x.$ Hence $e^{-x^{2}}v\left(  H_{\nu}\right)  $ is increasing
for all real $x.$ Thus
\[
e^{-x^{2}}v\left(  H_{\nu}\right)  \left(  x\right)  <\lim_{x\rightarrow
\infty}e^{-x^{2}}v\left(  H_{\nu}\right)  \left(  x\right)  \leq0
\]
and it follows that $v\left(  H_{\nu}\right)  \left(  x\right)  <0$ for all
real $x.$
\end{proof}

\section{Generalized Laguerre polynomials}

The generalized Laguerre polynomial $L_{n}^{\left(  \alpha\right)  }$ is a
solution of the differential equation
\[
f^{\prime\prime}+\frac{\alpha+1-x}{x}f^{\prime}+\frac{n}{x}f=0
\]
where $\alpha$ is a complex number. Note that by (\ref{eq1})
\begin{equation}
v\left(  L_{n}^{\left(  \alpha\right)  }\right)  =\left(  \frac{d}{dx}%
L_{n}^{\left(  \alpha\right)  }\right)  ^{2}+\frac{\alpha+1-x}{x}%
L_{n}^{\left(  \alpha\right)  }\frac{d}{dx}L_{n}^{\left(  \alpha\right)
}+\frac{n}{x}\left(  L_{n}^{\left(  \alpha\right)  }\right)  ^{2}.
\label{eqvLaguerre}%
\end{equation}
It is well known that $L_{n}^{\left(  \alpha\right)  }$ is an orthogonal
polynomial for $\alpha>-1$, so $v(L_{n}^{\left(  \alpha\right)  })$ is
positive on the real line for $\alpha>-1.$ Alternatively one can try to apply
Theorem \ref{Thm2Monotdiffeq} and in this case the condition (\ref{eqcritb})
is satisfied exactly when the same condition $\alpha>-1$ holds.

G. Szeg\"{o} proved in \cite{Szeg48} the Tur\'{a}n inequality for normalized
Laguerre polynomials $L_{n}^{\left(  \alpha\right)  }/L_{n}^{\left(
\alpha\right)  }\left(  0\right)  $ and $\alpha>-1$ in the framework of the
Laguerre-P\'{o}lya class. More elementary proofs have been given in
\cite{RaTh49}, \cite{Szas50} and \cite{ThNa51} based on recurrence relations
and identities. In order to relate the Tur\'{a}n inequality to the function
$v(L_{n}^{\left(  \alpha\right)  })$ let us recall the recurrence relations
\begin{align*}
x\frac{d}{dx}L_{n}^{\left(  \alpha\right)  }\left(  x\right)   &
=nL_{n}^{\left(  \alpha\right)  }\left(  x\right)  -\left(  \alpha+n\right)
L_{n-1}^{\left(  \alpha\right)  }\left(  x\right) \\
x\frac{d}{dx}L_{n}^{\left(  \alpha\right)  }\left(  x\right)   &  =-\left(
\alpha+n+1-x\right)  L_{n}^{\left(  \alpha\right)  }\left(  x\right)  +\left(
n+1\right)  L_{n+1}^{\left(  \alpha\right)  }\left(  x\right)  .
\end{align*}
Then the product $-\left(  \alpha+n\right)  \left(  n+1\right)  L_{n+1}%
^{\left(  \alpha\right)  }L_{n-1}^{\left(  \alpha\right)  }$ is equal to
\[
x^{2}(\frac{d}{dx}L_{n}^{\left(  \alpha\right)  })^{2}+\left(  \alpha
+1-x\right)  xL_{n}^{\left(  \alpha\right)  }\frac{d}{dx}L_{n}^{\left(
\alpha\right)  }-n\left(  \alpha+n+1-x\right)  L_{n}^{\left(  \alpha\right)
}L_{n}^{\left(  \alpha\right)  }.
\]
Adding $n\left(  n+\alpha+1\right)  L_{n}^{\left(  \alpha\right)  }%
L_{n}^{\left(  \alpha\right)  }$ to both sides leads to
\begin{align*}
&  n\left(  n+\alpha+1\right)  L_{n}^{\left(  \alpha\right)  }L_{n}^{\left(
\alpha\right)  }-\left(  \alpha+n\right)  \left(  n+1\right)  L_{n+1}^{\left(
\alpha\right)  }L_{n-1}^{\left(  \alpha\right)  }\\
&  =x^{2}\left(  \frac{d}{dx}L_{n}^{\left(  \alpha\right)  }\right)
^{2}+\left(  \alpha+1-x\right)  xL_{n}^{\left(  \alpha\right)  }\frac{d}%
{dx}L_{n}^{\left(  \alpha\right)  }+nxL_{n}^{\left(  \alpha\right)  }%
L_{n}^{\left(  \alpha\right)  }.
\end{align*}
In view of (\ref{eqvLaguerre}) we obtain
\[
x^{2}v(L_{n}^{\left(  \alpha\right)  })=n\left(  n+\alpha+1\right)
L_{n}^{\left(  \alpha\right)  }L_{n}^{\left(  \alpha\right)  }-\left(
\alpha+n\right)  \left(  n+1\right)  L_{n+1}^{\left(  \alpha\right)  }%
L_{n-1}^{\left(  \alpha\right)  }.
\]
For $\alpha=0$ we see that
\[
x^{2}v(L_{n}^{\left(  0\right)  })=n\left(  n+1\right)  \Delta(L_{n}^{\left(
0\right)  })
\]
which is already mentioned in \cite{MuNa51}, see also \cite[p. 71]{Skov54}.
For the normalized $\widetilde{L}_{n}^{\left(  \alpha\right)  }\left(
x\right)  =L_{n}^{\left(  \alpha\right)  }\left(  x\right)  /L_{n}^{\left(
\alpha\right)  }\left(  0\right)  $ one obtains the formula
\[
x^{2}v(L_{n}^{\left(  \alpha\right)  })=L_{n}^{\left(  \alpha\right)  }\left(
0\right)  L_{n}^{\left(  \alpha\right)  }\left(  0\right)  n\left(
\alpha+1+n\right)  \Delta(\widetilde{L}_{n}^{\left(  \alpha\right)  })
\]
using that $L_{n}^{\left(  \alpha\right)  }\left(  0\right)  =\frac{\left(
\alpha+1\right)  _{n}}{n!}$.

\section{Ultraspherical polynomials}

The ultraspherical polynomial of degree $n$ is denoted by $P_{n}^{\left(
\lambda\right)  }\left(  x\right)  ,$ and the normalized ultraspherical
polynomial is denoted by
\[
p_{n}^{\left(  \lambda\right)  }\left(  x\right)  =\frac{P_{n}^{\left(
\lambda\right)  }\left(  x\right)  }{P_{n}^{\left(  \lambda\right)  }\left(
1\right)  }.
\]
Then $p_{n}^{\left(  \lambda\right)  }$ and $P_{n}^{\left(  \lambda\right)  }$
satisfy the differential equation
\[
f^{\prime\prime}\left(  x\right)  -\frac{2\lambda+1}{1-x^{2}}xf^{\prime
}\left(  x\right)  +\frac{n\left(  n+2\lambda\right)  }{1-x^{2}}f=0.
\]
In this section we will give a new proof of the fact that
\[
\Delta(p_{n}^{\left(  \lambda\right)  })=p_{n}^{\left(  \lambda\right)  }%
p_{n}^{\left(  \lambda\right)  }-p_{n+1}^{\left(  \lambda\right)  }%
p_{n-1}^{\left(  \lambda\right)  }\geq0\text{ on }\left[  -1,1\right]  .
\]
We commence with the observation that $\Delta(p_{n}^{\left(  \lambda\right)
})$ is up to a positive factor equal to $v_{g}(p_{n}^{\left(  \lambda\right)
})$ for a suitable function $g.$

\begin{theorem}
\label{ThmUltraTuran}Define $g\left(  x\right)  =\left(  1-x^{2}\right)
^{\frac{1}{2}}$. Then the following identity holds:
\[
\Delta(p_{n}^{\left(  \lambda\right)  })=\frac{\left(  1-x^{2}\right)  ^{2}%
}{n\left(  n+2\lambda\right)  }v_{g}(p_{n}^{\left(  \lambda\right)  }).
\]

\end{theorem}

\begin{proof}
In the proof we write $p_{n}$ instead of $p_{n}^{\left(  \lambda\right)  }.$
Then (\ref{eq1}) shows that
\[
v\left(  p_{n}\right)  =p_{n}^{\prime}p_{n}^{\prime}-\frac{2\lambda+1}%
{1-x^{2}}xp_{n}^{\prime}p_{n}+\frac{n\left(  n+2\lambda\right)  }{1-x^{2}%
}p_{n}^{2}.
\]
Using the recurrence relations (cf. \cite[p. 345]{Niko23})
\begin{align*}
p_{n+1}\left(  x\right)   &  =xp_{n}\left(  x\right)  -\frac{1-x^{2}%
}{n+2\lambda}p_{n}^{\prime}\left(  x\right) \\
p_{n-1}\left(  x\right)   &  =xp_{n}\left(  x\right)  +\frac{1-x^{2}}{n}%
p_{n}^{\prime}\left(  x\right)
\end{align*}
it is easy to see
\begin{equation}
-p_{n+1}p_{n-1}=\frac{\left(  1-x^{2}\right)  ^{2}}{n\left(  n+2\lambda
\right)  }p_{n}^{\prime}p_{n}^{\prime}-\frac{2\lambda x\left(  1-x^{2}\right)
}{n\left(  n+2\lambda\right)  }p_{n}^{\prime}p_{n}-x^{2}p_{n}^{2}.
\label{eqhalfTuran}%
\end{equation}
Adding $p_{n}^{2}$ on both sides show that
\[
\Delta\left(  p_{n}\right)  =\frac{\left(  1-x^{2}\right)  ^{2}}{n\left(
n+2\lambda\right)  }p_{n}^{\prime}p_{n}^{\prime}-\frac{2\lambda x\left(
1-x^{2}\right)  }{n\left(  n+2\lambda\right)  }p_{n}^{\prime}p_{n}+\left(
1-x^{2}\right)  p_{n}^{2}.
\]
It follows that
\begin{align*}
\Delta\left(  p_{n}\right)   &  =\frac{\left(  1-x^{2}\right)  ^{2}}{n\left(
n+2\lambda\right)  }\left(  p_{n}^{\prime}p_{n}^{\prime}-\frac{2\lambda
x}{1-x^{2}}p_{n}^{\prime}p_{n}+\frac{n\left(  n+2\lambda\right)  }{1-x^{2}%
}p_{n}^{2}\right) \\
&  =\frac{\left(  1-x^{2}\right)  ^{2}}{n\left(  n+2\lambda\right)  }\left(
v\left(  p_{n}\right)  +\frac{x}{1-x^{2}}p_{n}^{\prime}p_{n}\right)  .
\end{align*}
Define $g\left(  x\right)  =\left(  1-x^{2}\right)  ^{\frac{1}{2}}$. Then
\[
\frac{g^{\prime}}{g}=\frac{\frac{d}{dx}\left(  1-x^{2}\right)  ^{\frac{1}{2}}%
}{\left(  1-x^{2}\right)  ^{\frac{1}{2}}}=\frac{-x}{1-x^{2}}%
\]
and the identity $v_{g}\left(  p_{n}\right)  =v\left(  p_{n}\right)
-\frac{g^{\prime}}{g}p_{n}^{\prime}p_{n}$ shows that
\[
\Delta\left(  p_{n}\right)  =\frac{\left(  1-x^{2}\right)  ^{2}}{n\left(
n+2\lambda\right)  }\left(  v\left(  p_{n}\right)  -\frac{g^{\prime}}{g}%
p_{n}^{\prime}p_{n}\right)  =\frac{\left(  1-x^{2}\right)  ^{2}}{n\left(
n+2\lambda\right)  }v_{g}\left(  p_{n}\right)  .
\]

\end{proof}

\begin{theorem}
\label{Thmvgppositive1}Let $g\left(  x\right)  =\left(  1-x^{2}\right)
^{\frac{1}{2}}$ and $q_{n}$ be a polynomial of degree $n$ with only real zeros
and contained in $\left[  -1,1\right]  $. Then $v_{g}\left(  q_{n}\right)
\left(  x\right)  \geq0$ for all $x\in\left(  -1,1\right)  .$
\end{theorem}

\begin{proof}
By Theorem \ref{ThmProduct} it suffices to show that $v_{g}\left(  f\right)
\left(  x\right)  \geq0$ for all $x\in\left[  -1,1\right]  $ where $f\left(
x\right)  =x-a$ and $a\in\left(  -1,1\right)  .$ Clearly
\[
v_{g}\left(  f\right)  \left(  x\right)  =-\left(  1-x^{2}\right)  ^{-\frac
{1}{2}}\left(  x-a\right)  ^{2}\frac{d}{dx}\frac{\left(  1-x^{2}\right)
^{\frac{1}{2}}\frac{d}{dx}\left(  x-a\right)  }{x-a}=\frac{1-ax}{1-x^{2}}.
\]
The numerator is a linear polynomial which is positive on the endpoints since
$1-a\geq0$ and $1+a\geq0.$
\end{proof}

Combining Theorem \ref{ThmUltraTuran} and Theorem \ref{Thmvgppositive1} we obtain:

\begin{theorem}
The Tur\'{a}n inequality holds for the ultraspherical polynomials:
\[
\Delta(p_{n}^{\left(  \lambda\right)  })\left(  x\right)  \geq0\text{ for all
}x\in\left[  -1,1\right]  .
\]

\end{theorem}

In the next section we will need the following extension of Theorem
\ref{Thmvgppositive1}:

\begin{proposition}
\label{Propvgpositive}Let $h\in C\left[  -1,1\right]  $ and $H=\int
\frac{h\left(  x\right)  }{1-x^{2}}$ and $g=e^{H}.$ Then the function
$v_{g}\left(  p\right)  $ is non-negative on $\left(  -1,1\right)  $ for all
odd and even polynomials with only real zeros, contained in $\left[
-1,1\right]  $ if and only if
\begin{equation}
1-xh\left(  x\right)  -x^{2}\geq0\text{ and }1+xh\left(  x\right)  +x^{2}%
\geq0. \label{eqsuffcond}%
\end{equation}

\end{proposition}

\begin{proof}
Since $g=e^{H}$ we have
\[
\frac{g^{\prime}\left(  x\right)  }{g\left(  x\right)  }=H^{\prime}\left(
x\right)  =\frac{h\left(  x\right)  }{1-x^{2}}.
\]
In the following we use the identity
\begin{equation}
v_{g}\left(  f\right)  \left(  x\right)  =v\left(  f\right)  -\frac{g^{\prime
}}{g}f^{\prime}f. \label{eqvgnew}%
\end{equation}
Define $e_{1}\left(  x\right)  =x.$ Clearly $v\left(  e_{1}\right)  =1$ and by
(\ref{eqvgnew})
\begin{equation}
v_{g}\left(  e_{1}\right)  \left(  x\right)  =1-\frac{h\left(  x\right)
}{1-x^{2}}x=\frac{1-xh\left(  x\right)  -x^{2}}{1-x^{2}}. \label{eqvg1}%
\end{equation}
For $f_{a}\left(  x\right)  =\left(  x-a\right)  \left(  x+a\right)  $ a
simple calculation shows that $v\left(  f_{a}\right)  =2a^{2}+2x^{2}$ and by
(\ref{eqvgnew})
\begin{equation}
v_{g}\left(  f_{a}\right)  \left(  x\right)  =2a^{2}+2x^{2}-\frac{h\left(
x\right)  }{1-x^{2}}2x\left(  x^{2}-a^{2}\right)  . \label{eqvg2}%
\end{equation}
For the necessity observe that $v_{g}\left(  f\right)  \left(  x\right)
\geq0$ for $f=e_{1}$ implies that $1-xh\left(  x\right)  -x^{2}\geq0$ by
(\ref{eqvg1}), and $v_{g}\left(  f\right)  \left(  x\right)  \geq0$ for
$f=f_{\alpha}$ and $a=1$ implies by (\ref{eqvg2}) that
\[
v_{g}\left(  f_{1}\right)  =2+2x^{2}+2xh\left(  x\right)  \geq0.
\]
For the converse, note that by Theorem \ref{ThmProduct} it suffices to show
that $v_{g}\left(  f\right)  \left(  x\right)  \geq0$ for all $x\in\left[
-1,1\right]  $ where $f=e_{1}$ and $f=f_{\alpha}.$ By (\ref{eqvg1}) and our
assumption ~$1-xh\left(  x\right)  -x^{2}\geq0$ we infer that $v_{g}\left(
e_{1}\right)  \geq0$ on $\left[  -1,1\right]  .$ For $f_{a}$ note that
\[
\left(  1-x^{2}\right)  v_{g}\left(  f_{a}\right)  =2\left(  1+xh\left(
x\right)  -x^{2}\right)  a^{2}+2x^{2}\left(  1-xh\left(  x\right)
-x^{2}\right)  .
\]
By assumption we know that $1-xh\left(  x\right)  -x^{2}\geq0.$When
$1+xh\left(  x\right)  -x^{2}\geq0$ it follows that the sum is positive, so
$v_{g}\left(  f_{a}\right)  \geq0.$ If $1+xh\left(  x\right)  -x^{2}$ is
negative, then the expression is concave function in the variable $a\in\left[
-1,1\right]  ,$ so we have to check the positivity only in the endpoints, i.e.
that we have to show that
\[
2\left(  1+xh\left(  x\right)  -x^{2}\right)  +2x^{2}\left(  1-xh\left(
x\right)  -x^{2}\right)  =2\left(  1-x^{2}\right)  \left(  x^{2}+xh\left(
x\right)  +1\right)
\]
is non-negative. This follows from our second assumption that $x^{2}+xh\left(
x\right)  +1\geq0.$
\end{proof}

In the last theorem we have required the condition $v_{g}\left(  p_{n}\right)
$ for all even or odd polynomials with only real zeros. If we require the
inequality for all polynomials with real zeros we obtain sharper requirements:

\begin{theorem}
Let $h\in C\left[  -1,1\right]  $ and $H=\int\frac{h\left(  x\right)
}{1-x^{2}} $ and $g=e^{H}.$ Then the function $v_{g}\left(  p\right)  $ is
non-negative on $\left(  -1,1\right)  $ for all polynomials with only real
zeros, contained in $\left[  -1,1\right]  $ if and only if
\[
1+x+h\left(  x\right)  \geq0\text{ and }1-x-h\left(  x\right)  \geq0
\]
for all $x\in\left[  -1,1\right]  .$
\end{theorem}

\begin{proof}
Let $f_{a}\left(  x\right)  =x-a$ with $a\in\left[  -1,1\right]  .$ Then by
(\ref{eqvgnew})
\[
v_{g}\left(  f_{a}\right)  \left(  x\right)  =1-\frac{h\left(  x\right)
}{1-x^{2}}\left(  x-a\right)  =\frac{1-x^{2}+\left(  a-x\right)  h\left(
x\right)  }{1-x^{2}}.
\]
For the necessity note that $v_{g}\left(  f_{a}\right)  \geq0$ on $\left[
-1,1\right]  $ for $a=1$ and $a=-1$ implies that
\begin{align*}
1-x^{2}+\left(  1-x\right)  h\left(  x\right)   &  \geq0\\
1-x^{2}-\left(  1+x\right)  h\left(  x\right)   &  \geq0.
\end{align*}
Division by $1-x$ of the first equation and division by $1+x$ of the second
equation give the necessary conditions.

For the converse, it suffices by Theorem \ref{ThmProduct} to show that
$v_{g}\left(  f_{a}\right)  \left(  x\right)  \geq0$ for all $x\in\left[
-1,1\right]  $ and all $a\in\left[  -1,1\right]  .$ Note that
\[
\left(  1-x^{2}\right)  v_{g}\left(  f\right)  =1-x^{2}+\left(  a-x\right)
h\left(  x\right)  ,
\]
and in the variable $a$ this is a linear polynomial. It is non-negative on
$\left[  -1,1\right]  $ if it is non-negative for $a=1$ and $a=-1,$ which are
equivalent to our necessary conditions.
\end{proof}

Let us finish this section with some comments. In \cite[p. 381]{ThNa51} the
authors showed that for the normalized ultraspherical polynomial
$p_{n}^{\left(  \lambda\right)  }$ and for $n>1$ the function
\[
\varphi_{n,\lambda}\left(  x\right)  =\frac{\Delta(p_{n}^{\left(
\lambda\right)  })\left(  x\right)  }{1-x^{2}}%
\]
increases on $\left[  0,1\right]  $ for $\lambda>0.$ Thus, using that
$\varphi_{n,\lambda}$ is even, one has $\varphi_{n,\lambda}\left(  0\right)
\leq\varphi_{n,\lambda}\left(  x\right)  \leq\varphi_{n,\lambda}\left(
1\right)  $ for all $x\in\mathbb{R}$. It follows (as pointed out in \cite[p.
381]{ThNa51}) that
\begin{equation}
\Delta(p_{n}^{\left(  \lambda\right)  })\left(  0\right)  \cdot\left(
1-x^{2}\right)  \leq\Delta(p_{n}^{\left(  \lambda\right)  })\left(  x\right)
\leq\frac{1-x^{2}}{1+2\lambda} \label{ineqThNa}%
\end{equation}
holds where the rule of L'Hospital is used for the determination of the
constant
\[
\lim_{x\rightarrow1}\frac{\Delta(p_{n}^{\left(  \lambda\right)  })\left(
x\right)  }{1-x^{2}}=-\frac{1}{2}\frac{d}{dx}\Delta(p_{n}^{\left(
\lambda\right)  })\left(  1\right)  =\frac{1}{1+2\lambda}.
\]
The inequality (\ref{ineqThNa}) has been rediscovered in \cite{Cons05} and
\cite{AGKL07} for the case of Legendre polynomials. For $\lambda\in(-\frac
{1}{2},0)$ the function $\varphi_{n,\lambda}$ decreases on $\left[
0,1\right]  $ and the inequality is reversed:
\[
\frac{1-x^{2}}{1+2\lambda}\leq\Delta\left(  p_{n}^{\left(  \lambda\right)
}\right)  \left(  x\right)  \leq\Delta\left(  p_{n}^{\left(  \lambda\right)
}\right)  \left(  0\right)  \cdot\left(  1-x^{2}\right)  .
\]
In passing we mention that G. Nikolov proved in \cite{Niko23} that
$\varphi_{n,\lambda}$ is concave for $\lambda\in\left(  -\frac{1}{2},0\right)
$ and convex for $\lambda>0.$

\section{The Tur\'{a}n inequality of Krasikov}

In the following we want to prove a Tur\'{a}n inequality due to I. Krasikov
\cite{Kras23}:
\begin{equation}
\Delta_{c\left(  \lambda\right)  }(p_{n}^{\left(  \lambda\right)  })\left(
x\right)  :=\left\vert x\right\vert ^{c\left(  \lambda\right)  }\left(
p_{n}^{\left(  \lambda\right)  }\left(  x\right)  \right)  ^{2}-p_{n+1}%
^{\left(  \lambda\right)  }\left(  x\right)  p_{n-1}^{\left(  \lambda\right)
}\left(  x\right)  \geq0 \label{eqDeltac}%
\end{equation}
for all $x\in\left[  -1,1\right]  ,$ and $\lambda>-\frac{1}{2}$ where the
exponent $c\left(  \lambda\right)  $ is defined by
\[
c\left(  \lambda\right)  =\frac{2}{1+2\left\vert \lambda\right\vert }\leq2.
\]
We shall offer here a shorter and simpler proof. Let us emphasize that for
$\lambda\in\left(  -\frac{1}{2},0\right)  ,$ Krasikov obtained a better bound
and it seems that our methods are limited. Clearly the estimate $c\left(
\lambda\right)  \geq1$ holds for all $\left\vert \lambda\right\vert \leq
\frac{1}{2}$ and Krasikov's theorem then implies the estimate of G. Nikolov
and V. Pillwein \cite{NiPi15}
\[
\left\vert x\right\vert \left(  p_{n}^{\left(  \lambda\right)  }\right)
^{2}-p_{n+1}^{\left(  \lambda\right)  }\left(  x\right)  p_{n-1}^{\left(
\lambda\right)  }\left(  x\right)  \geq0
\]
for all $\left\vert \lambda\right\vert \leq\frac{1}{2}$ and $x\in\left[
-1,1\right]  .$ This inequality was earlier proved by Gerhold and Kauers in
\cite{GeKa06} for the case of Legendre polynomials, so for $\lambda=\frac
{1}{2}.$

For $x\in\left[  0,1\right]  $ and $\alpha>0$ the inequality $0\leq x^{\alpha
}\leq1$ is obvious. It follows that for $c\in\left[  0,2\right]  $ and
$x\in\left[  0,1\right]  $
\[
x^{c}-x^{2}=x^{c}\left(  1-x^{2-c}\right)  \in\left[  0,1\right]  \text{.}%
\]
We now proceed to the proof that $\Delta_{c}\left(  p_{n}\right)  =\left\vert
x\right\vert ^{c}p_{n}^{2}-p_{n+1}p_{n-1}$ is non-negative on $\left[
-1,1\right]  $ for a suitable $c\in\left(  0,2\right)  .$ Add in
(\ref{eqhalfTuran}) the term $\left\vert x\right\vert ^{c}p_{n}^{2}$ on both
sides, then
\[
\Delta_{c}\left(  p_{n}\right)  =\frac{\left(  1-x^{2}\right)  ^{2}}{n\left(
n+2\lambda\right)  }p_{n}^{\prime}p_{n}^{\prime}-\frac{2\lambda x\left(
1-x^{2}\right)  }{n\left(  n+2\lambda\right)  }p_{n}^{\prime}p_{n}+\left(
\left\vert x\right\vert ^{c}-x^{2}\right)  p_{n}^{2}.
\]
Then
\[
\Delta_{c}\left(  p_{n}\right)  =\frac{\left(  1-x^{2}\right)  \left(
\left\vert x\right\vert ^{c}-x^{2}\right)  }{n\left(  n+2\lambda\right)
}D_{c}\left(  p_{n}\right)
\]
where
\begin{align*}
D_{c}\left(  p_{n}\right)   &  =\frac{1-x^{2}}{\left\vert x\right\vert
^{c}-x^{2}}p_{n}^{\prime}p_{n}^{\prime}-\frac{2\lambda x\left(  1-x^{2}%
\right)  }{\left(  \left\vert x\right\vert ^{c}-x^{2}\right)  }\frac
{p_{n}^{\prime}p_{n}}{1-x^{2}}+\frac{n\left(  n+2\lambda\right)  }{1-x^{2}%
}p_{n}^{2}\\
&  \geq p_{n}^{\prime}p_{n}^{\prime}-\frac{2\lambda xp_{n}^{\prime}p_{n}%
}{\left\vert x\right\vert ^{c}-x^{2}}+\frac{n\left(  n+2\lambda\right)
}{1-x^{2}}p_{n}^{2}%
\end{align*}
using the fact that for all $x\in\left[  0,1\right]  $%
\[
\frac{1-x^{2}}{\left\vert x\right\vert ^{c}-x^{2}}\geq\frac{2}{2-c}\geq1.
\]
The latter inequality holds since the left hand side is a decreasing function
and the right hand side is equal to the evaluation at $x=1.$ Now define%
\begin{equation}
h_{\lambda,c}\left(  x\right)  =\frac{2\lambda x\left(  1-x^{2}\right)
}{\left\vert x\right\vert ^{c}-x^{2}}-\left(  2\lambda+1\right)  x.
\label{eqdefhlc}%
\end{equation}
It follows
\begin{align*}
D_{c}\left(  p_{n}\right)   &  \geq p_{n}^{\prime}p_{n}^{\prime}-\frac{\left(
2\lambda+1\right)  p_{n}^{\prime}p_{n}}{1-x^{2}}-\frac{h_{\lambda,c}\left(
x\right)  }{1-x^{2}}p_{n}^{\prime}p_{n}+\frac{n\left(  n+2\lambda\right)
}{1-x^{2}}p_{n}^{2}\\
&  =v\left(  p_{n}\right)  -\frac{h_{\lambda,c}\left(  x\right)  }{1-x^{2}%
}p_{n}^{\prime}p_{n}=v_{g}\left(  p_{n}\right)
\end{align*}
where in the last step we define $g_{\lambda,c}=e^{H_{\lambda,c}}$ and
$H_{\lambda,c}=\int\frac{h_{\lambda,c}\left(  x\right)  }{1-x^{2}}dx$, so
$\frac{g^{\prime}}{g}=\frac{h_{\lambda,c}\left(  x\right)  }{1-x^{2}}.$

Now we are going to prove the Tur\'{a}n inequality of Krasikov:

\begin{theorem}
Let $\lambda>-\frac{1}{2}$ and $c\left(  \lambda\right)  =\frac{2}%
{1+2\left\vert \lambda\right\vert }$. Then the Tur\'{a}n inequality
(\ref{eqDeltac}) holds.
\end{theorem}

\begin{proof}
By the above arguments it is suffices to show that for $c=c\left(
\lambda\right)  $
\[
v_{g_{\lambda,c}}\left(  p_{n}\right)  \left(  x\right)  \geq0\text{ for all
}x\in\left[  -1,1\right]
\]
for each ultraspherical polynomial $p_{n}=p_{n}^{\left(  \lambda\right)  }$.
Here $g_{\lambda,c}=e^{H_{\lambda,c}}$ and $H_{\lambda,c}=\int\frac
{h_{\lambda,c}\left(  x\right)  }{1-x^{2}}dx.$ By Proposition
\ref{Propvgpositive} it suffices to prove that%
\[
1-xh_{\lambda,c}\left(  x\right)  -x^{2}\geq0\text{ and }1+xh_{\lambda
,c}\left(  x\right)  +x^{2}\geq0.
\]
Define $H_{\delta}^{\lambda}\left(  x\right)  =1+x\delta h_{\lambda,c}\left(
x\right)  +\delta x^{2}.$ Then it suffices to show that $H_{\delta}^{\lambda
}\left(  x\right)  \geq0$ for all $x\in\left[  -1,1\right]  $ and $\delta
=\pm1.$ We can compute $H_{\delta}^{\lambda}$ by using the definition of
$h_{\lambda,c}$ in (\ref{eqdefhlc}) and we obtain
\begin{align*}
H_{\delta}^{\lambda}\left(  x\right)   &  =1+\frac{2\lambda\delta x^{2}\left(
1-x^{2}\right)  }{\left\vert x\right\vert ^{c}-x^{2}}-\left(  2\lambda
+1\right)  \delta x^{2}+\delta x^{2}\\
&  =\frac{\left\vert x\right\vert ^{c}-x^{2}+2\lambda\delta x^{2}\left(
1-x^{2}\right)  -2\lambda x^{2}\delta\left(  \left\vert x\right\vert
^{c}-x^{2}\right)  }{\left\vert x\right\vert ^{c}-x^{2}}\\
&  =\frac{\left\vert x\right\vert ^{c}-x^{2}+2\lambda\delta x^{2}%
-2\lambda\delta x^{2}\left\vert x\right\vert ^{c}}{\left\vert x\right\vert
^{c}-x^{2}}.
\end{align*}
The last formula shows that $H_{\delta}^{\lambda}\left(  x\right)  $ is an
even function. Hence it suffices to show that
\[
R_{\delta}^{\lambda}\left(  x\right)  :=x^{c}-x^{2}+2\lambda\delta
x^{2}-2\lambda\delta x^{2+c}=x^{c}\left(  1-x^{2-c}\right)  +2\lambda\delta
x^{2}\left(  1-x^{c}\right)  .
\]
is non-negative on $\left[  0,1\right]  $ for $\delta=-1,1.$ Let us write
$R_{\delta}^{\lambda}\left(  x\right)  =x^{c}r_{\delta}^{\lambda}\left(
x\right)  $ with
\[
r_{\delta}^{\lambda}\left(  x\right)  =1-x^{2-c}+2\lambda\delta x^{2-c}%
-2\lambda\delta x^{2}.
\]
Note that $r_{\delta}^{\lambda}\left(  0\right)  =1$ and $r_{\delta}^{\lambda
}\left(  1\right)  =0.$In the case $\lambda=0,$ the function $r_{\delta}^{0}$
is obviously non-negative for all $x\in\left[  0,1\right]  .$ For
$\lambda\delta\geq0$ it is clear that $r_{\delta}^{\lambda}\left(  x\right)
\geq0$ for all $x\in\left[  0,1\right]  $ since $x^{2-c}$ and $x^{c}\in\left[
0,1\right]  .$ Assume that $\lambda\delta\leq0.$ Then
\[
\frac{d}{dx}r_{\delta}^{\lambda}\left(  x\right)  =\left(  2\lambda
\delta-1\right)  \left(  2-c\right)  x^{1-c}-4\lambda\delta x.
\]
It follows that $r_{\delta}^{\lambda}$ has exactly one critical point at
$x_{0}>0$ where $x_{0}$ is given by
\[
x_{0}^{c}=\frac{\left(  2\lambda\delta-1\right)  \left(  2-c\right)
}{4\lambda\delta}.
\]
In the case that
\[
c\leq\frac{2}{1-2\lambda\delta}=\frac{2}{1+2\left\vert \lambda\right\vert }%
\]
one sees that
\[
x_{0}^{c}=\frac{\left(  2\lambda\delta-1\right)  \left(  2-c\right)
}{4\lambda\delta}\geq\frac{\left(  2\lambda\delta-1\right)  \left(  2-\frac
{2}{1-2\lambda\delta}\right)  }{4\lambda\delta}=1.
\]
Thus $x_{0}\geq1,$ and it follows that $r_{\delta}^{\lambda}\left(  x\right)
$ is monotone increasing or decreasing on $\left[  0,1\right]  .$ Since
$r_{\delta}^{\lambda}\left(  0\right)  =1$ and $r_{\delta}^{\lambda}\left(
1\right)  =0$ the function must be decreasing, and since $r_{\delta}^{\lambda
}\left(  1\right)  =0,$ it follows that $r_{\delta}^{\lambda}\left(  x\right)
\geq0$ on $\left[  0,1\right]  .$
\end{proof}

Finally let us mention that Tur\'{a}n inequalities for general orthogonal
polynomials, based on an analysis of the recurrence relation, have been
investigated by many authors, and we just mention  some recent  literature in
\cite{Kras05}, \cite{Kras11}, \cite{Szwa21}.

\textbf{Acknowledgement} The work of all authors was funded under projects
BG-175467353-2024-18-0018-C01  and  KP-06-RILA  (entry number BG-175467353-2024-07-0008) 
with Bulgarian NSF.

\smallskip\noindent\textbf{Data availability} Data sharing is not applicable
to this article since no data sets were generated or analyzed.

\smallskip\noindent\textbf{Declarations:}

\smallskip\noindent\textbf{Conflict of interest} The authors declare that they
have no competing interests.

\smallskip\noindent\textbf{Usage declaration}: We have not used generative AI
.

\end{document}